\documentclass[a4paper,twoside,11pt]{article}
\usepackage{bbm}
\usepackage{amsfonts}
\usepackage{amssymb}\usepackage{amsmath}
 \numberwithin{equation}{section}
\begin{document}

\title{Blow-up phenomena and global existence for the periodic two-component
Dullin-Gottwald-Holm system}

\author{Jingjing Liu\footnote{E-mail: jingjing830306@163.com.}
\\Department of Mathematics and Information Science, \\ Zhengzhou
University of Light Industry , 450002 Zhengzhou, China}
\date{}
\maketitle

\begin{abstract}
This paper is concerned with blow-up phenomena and global existence
for the periodic two-component Dullin-Gottwald-Holm system. We first
obtain several blow-up results and the blow-up rate of strong
solutions to the system. We then present a global existence result
for strong solutions to the system.\\


\end{abstract}

\section{Introduction}
\noindent

In this paper, we consider the following periodic two-component
Dullin-Gottwald-Holm (DGH) system:
\begin{equation}\label{eq:original}
\left\{\begin{array}{ll}m_{t}-Au_{x}+um_{x}+2u_{x}m+\gamma u_{xxx}+\rho\rho_{x}=0,& t > 0,\,x\in \mathbb{R},\\
 \rho_{t}+(u\rho )_{x}=0, &t > 0,\,x\in \mathbb{R},\\
u(0,x)=u_{0}(x),&x\in \mathbb{R}, \\
\rho(0,x) = \rho_{0}(x),&x\in \mathbb{R},\\
u(t,x+1)=u(t,x), &t \geq 0,\,x\in \mathbb{R},\\
\rho(t,x+1)=\rho(t,x), &t \geq 0,\,x\in \mathbb{R},\\
\end{array}\right.
\end{equation}
where $m=u-u_{xx},$ $A>0$ and $\gamma$ are constants.
\newline

The system (1.1) has been recently derived by Zhu et al. in
\cite{zhu} by follow Ivanov's approach \cite{Ivanov}. It was shown
in \cite{zhu} that the DGH system is completely integrable and can
be written as a compatibility condition of two linear systems
$$\Psi_{xx}=\left(-\xi^{2}\rho^{2}+\xi\left(m-\frac{A}{2}+\frac{\gamma}{2}\right)+\frac{1}{4}\right)\Psi$$
and
$$\Psi_{t}=\left(\frac{1}{2\xi}-u+\gamma\right)\Psi_{x}+\frac{1}{2}u_{x}\Psi,$$
where $\xi$ is a spectral parameter. Moreover, this system has the
following two Hamiltonians
$$E(u,\rho)=\frac{1}{2}\int(u^{2}+u_{x}^{2}+(\rho-1)^{2})dx$$
and
$$F(u,\rho)=\frac{1}{2}\int(u^{3}+uu_{x}^{2}-Au^{2}-\gamma
u_{x}^{2}+2u(\rho-1)+u(\rho-1)^{2})dx.$$\

For $\rho=0$ and $m=u-\alpha^{2}u_{xx},$ (1.1) becomes to the DGH
equation \cite{DGH}
$$
u_{t}-\alpha^{2}u_{txx}-Au_{x}+3uu_{x}+\gamma
u_{xxx}=\alpha^{2}(2u_{x}u_{xx}+uu_{xxx}),
$$
where $A$ and $\alpha$ are two positive constants, modeling
unidirectional propagation of surface waves on a shallow layer of
water which is at rest at infinity, $u(t,x)$ standing for fluid
velocity. It is completely integrable with a bi-Hamiltonian and a
Lax pair. Moreover, its traveling wave solutions include both the
KdV solitons and the CH peakons as limiting cases \cite{DGH}. The
Cauchy problem of the DGH equation has been extensively studied, cf.
\cite{yin,shen,ai,guo,meng,sun,liu,tian,yin2,zhang,zhang1,m,yan}.
\newline

For $\rho\not \equiv 0$, $\gamma=0$, the system (1.1) becomes to the
two-component Camassa-Holm system \cite{Ivanov}
\begin{equation}
\left\{\begin{array}{ll}m_{t}-Au_{x}+um_{x}+2u_{x}m+\rho\rho_{x}=0,\\
 \rho_{t}+(u\rho )_{x}=0,\\
\end{array}\right.
\end{equation}
where $\rho(t,x)$ in connection with the free surface elevation from
scalar density (or equilibrium) and the parameter $A$ characterizes
a linear underlying shear flow. The system (1.2) describes water
waves in the shallow water regime with nonzero constant vorticity,
where the nonzero vorticity case indicates the presence of an
underlying current. A large amount of literature was devoted to the
Cauchy problem (1.2), see \cite{chen1, chen2, gui1,gui2,
con,escher,fu,zhang2,hu}.
\newline

The Cauchy problem (1.1) has been discussed in \cite{zhu}. Therein
Zhu and Xu established the local well-posedness to the system (1.1),
derived the precise blow-up scenario and investigated the wave
breaking for the system (1.1). The aim of this paper is to study
further the blow-up phenomena for strong solutions to (1.1) and to
present a global existence result.
\newline

Our paper is organized as follows. In Section 2, we briefly give
some needed results including the local well posedness of the system
(1.1), the precise blow-up scenarios and some useful lemmas to study
blow-up phenomena and global existence. In Section 3, we give
several new blow-up results and the precise blow-up rate. In Section
4, we present a new global existence result of strong solutions
to (1.1).\\
\newline
\textbf{Notation}  Given a Banach space $Z$, we denote its norm by
 $\|\cdot\|_{Z}$. Since all space of functions are over
 $\mathbb{S}$, for simplicity, we drop $\mathbb{S}$ in our notations
 if there is no ambiguity.

\section{Preliminaries}
\noindent

In this section, we will briefly give some needed results in order
to pursue our goal. \

With $m=u-u_{xx}$, we can rewrite the system (1.1) as follows:
\begin{equation}
\left\{\begin{array}{ll}u_{t}-u_{txx}-Au_{x}+\gamma u_{xxx}+3uu_{x}-2u_{x}u_{xx}-uu_{xxx}+\rho\rho_{x}=0,& t > 0,\,x\in \mathbb{R},\\
 \rho_{t}+(u\rho )_{x}=0, &t > 0,\,x\in \mathbb{R},\\
u(0,x)=u_{0}(x),&x\in \mathbb{R}, \\
\rho(0,x) = \rho_{0}(x),&x\in \mathbb{R},\\
u(t,x+1)=u(t,x), &t \geq 0,\,x\in \mathbb{R},\\
\rho(t,x+1)=\rho(t,x), &t \geq 0,\,x\in \mathbb{R}.\\
\end{array}\right.
\end{equation}
Note that if $G(x):=\frac{\cosh(x-[x]-1/2)}{2\sinh(1/2)}$, $x\in
\mathbb{R}$ is the kernel of $(1-
\partial^{2}_{x})^{-1}$, then $(1-
\partial^{2}_{x})^{-1}f = G*f $ for all $f \in L^{2}(\mathbb{S})$, $G \ast
m=u $. Here we denote by $\ast$ the convolution. Using this
identity, we can rewrite the system (2.1) as follows:
\begin{equation}
\left\{\begin{array}{ll}u_{t}+(u-\gamma)u_{x}=
-\partial_{x}G*\left(u^{2}+\frac{1}{2}u_{x}^{2}+(\gamma-A)u+\frac{1}{2}\rho^{2}\right),& t > 0,\,x\in \mathbb{R},\\
\rho_{t}+(u\rho )_{x}=0, &t > 0,\,x\in \mathbb{R},\\
u(0,x)=u_{0}(x),&x\in \mathbb{R}, \\
\rho(0,x) = \rho_{0}(x),&x\in \mathbb{R},\\
u(t,x+1)=u(t,x), &t \geq 0,\,x\in \mathbb{R},\\
\rho(t,x+1)=\rho(t,x), &t \geq 0,\,x\in \mathbb{R},\\
\end{array}\right.
\end{equation}\

The local well-posedness of the Cauchy problem (2.1) can be obtained
by applying the Kato's theorem. As a result, we have the following
well-posedness result.\\
\newline
\textbf{Lemma 2.1.} \emph{(\cite{zhu}).} Given an initial data
$(u_{0},\rho_{0})\in H^{s}\times H^{s-1}, s\geq 2,$ there exists a
maximal $T = T (\|(u_{0},\rho_{0})\|_{H^{s}\times H^{s-1}})> 0$ and
a unique solution $$ (u,\rho)\in C([0,T); H^{s}\times H^{s-1})\cap
C^{1}([0,T);H^{s-1}\times H^{s-2})
$$ of (2.1).
Moreover, the solution $(u,\rho)$ depends continuously on the
initial data $(u_{0},\rho_{0})$ and the maximal time of existence $T
> 0$ is independent of $s$.\\
\

Consider now the following initial value problem
\begin{equation}
\left\{\begin{array}{ll}q_{t}=u(t,q),\ \ \ \ t\in[0,T), \\
q(0,x)=x,\ \ \ \ x\in\mathbb{R}, \end{array}\right.
\end{equation}
where $u$ denotes the first component of the solution $(u,\rho)$ to (2.1).\\
\newline
\textbf{Lemma 2.2.} \emph{(\cite{zhu}).} Let $(u, \rho)$ be the
solution of (2.1) with initial data $(u_{0},\rho_{0})\in H^{s}\times
H^{s-1}, s\geq 2$. Then Eq.(2.3) has a unique solution $q\in
C^1([0,T)\times \mathbb{R};\mathbb{R})$. Moreover, the map
$q(t,\cdot)$ is an increasing diffeomorphism of $\mathbb{R}$ with
$$
q_{x}(t,x)=exp\left(\int_{0}^{t}u_{x}(s,q(s,x))ds\right)>0, \ \
(t,x)\in [0,T)\times \mathbb{R}.$$
\newline
\textbf{Lemma 2.3.} \emph{(\cite{zhu}).} Let $(u, \rho)$ be the
solution of (2.1) with initial data $(u_{0},\rho_{0})\in H^{s}\times
H^{s-1}, s\geq 2,$ and $T>0$ be the maximal existence. Then we have
$$
\rho(t,q(t,x))q_{x}(t,x)=\rho_{0}(x), \ (t,x)\in [0,T)\times
\mathbb{S}.
$$
Moreover, if there exists a $x_{0}\in \mathbb{S}$ such that
$\rho_{0}(x_{0})=0$, then $\rho(t,q(t,x_{0}))=0$ for all
$t\in[0,T).$\\
\

Next, we will give two useful conservation laws of strong
solutions to (2.1).\\
\newline
\textbf{Lemma 2.4.} \emph{(\cite{zhu}).} Let $(u, \rho)$ be the
solution of (2.1) with initial data $(u_{0},\rho_{0})\in H^{s}\times
H^{s-1}, s\geq 2,$ and $T>0$ be the maximal existence. Then for all
$t\in [0, T ),$ we have
$$
\int_{\mathbb{S}} (u^2+u_x^2+\rho^2)dx=\int_{\mathbb{S}}
(u_0^2+u_{0,x}^2+\rho_0^2)dx:=E_{0}.
$$
\newline
\textbf{Lemma 2.5.} Let $(u, \rho)$ be the solution of (2.1) with
initial data $(u_{0},\rho_{0})\in H^{s}\times H^{s-1}, s\geq 2,$ and
$T>0$ be the maximal existence. Then for all $t\in [0, T ),$ we have
$$\int_{\mathbb{S}}u(t,x)dx=\int_{\mathbb{S}}u_{0}(x)dx.$$
\newline
\textbf{Proof.} By the first equation in (2.1), we have
\begin{eqnarray*}
\frac{d}{dt}\int_{\mathbb{S}}u(t,x)dx&=&\int_{\mathbb{S}}u_{t}dx\\
&=&\int_{\mathbb{S}}(u_{txx}+Au_{x}-\gamma
u_{xxx}-3uu_{x}+2u_{x}u_{xx}+uu_{xxx}-\rho\rho_{x})dx=0
\end{eqnarray*}
This completes the proof of the lemma. $\hfill{} \Box$ \\
\

Then, we state the following precise blow-up mechanism of (2.1).\\
\newline
\textbf{Lemma 2.6.} \emph{(\cite{zhu}).} Let $(u, \rho)$ be the
solution of (2.1) with initial data $(u_{0},\rho_{0})\in H^{s}\times
H^{s-1}, s\geq 2,$ and $T>0$ be the maximal existence. Then the
solution blows up in finite time if and only if
$$\liminf\limits_{t\rightarrow T^{-}}\{\inf\limits_{x\in\mathbb{S}}u_{x}(t,x)\}=-\infty.$$
\newline
\textbf{Lemma 2.7.} \emph{(\cite{Constantin 4}).} Let $t_{0}>0$ and
$v\in C^{1}([0,t_{0}); H^{2}(\mathbb{R}))$. Then for every
$t\in[0,t_{0})$ there exists at least one point $\xi(t)\in
\mathbb{R}$ with
$$ m(t):=\inf_{x\in \mathbb{R}}\{v_{x}(t,x)\}=v_{x}(t,\xi(t)),$$ and
the function $m$ is almost everywhere differentiable on $(0,t_{0})$
with $$ \frac{d}{dt}m(t)=v_{tx}(t,\xi(t)) \ \ \ \ a.e.\ on \
(0,t_{0}).$$
\newline
\textbf{Lemma 2.8.} \emph{(\cite{Y2}).} (i) For every $f \in
H^{1}(\mathbb{S})$, we have
$$\max_{x \in [0,1]}f^{2}(x) \le
\frac{e+1}{2(e-1)}\| f\|^{2}_{H^{1}},
$$
where the constant $\frac{e+1}{2(e-1)}$ is sharp.
\par(ii)
For every $f \in H^{3}(\mathbb{S})$, we have
$$\max_{x \in [0,1]}f^{2}(x) \le
c\| f\|^{2}_{H^{1}},
$$
with the best possible constant $c$ lying within the range
$(1,\frac{13}{12}]$. Moreover, the best constant $c$ is
$\frac{e+1}{2(e-1)}$.\\
\newline
\textbf{Lemma 2.9.} \emph{(\cite{hu1}).} If $f\in H^{3}(\mathbb{S})$
is such that $\int_{\mathbb{S}}f(x)dx=\frac{a_{0}}{2},$ then for
every $\epsilon >0,$ we have
$$\max_{x \in [0,1]}f^{2}(x) \leq
\frac{\epsilon+2}{24}\int_{\mathbb{S}}f_{x}^{2}dx+\frac{\epsilon+2}{4\epsilon}a_{0}^{2}.$$
Moreover,
$$\max_{x \in [0,1]}f^{2}(x) \leq
\frac{\epsilon+2}{24}\|f\|_{H^{1}(\mathbb{S})}^{2}+\frac{\epsilon+2}{4\epsilon}a_{0}^{2}.$$
\newline
\textbf{Lemma 2.10.} \emph{(\cite{zhou}).} Assume that a
differentiable function $y(t)$ satisfies
\begin{equation}
y^{\prime}(t)\leq -Cy^{2}(t)+K
\end{equation} with constants $C,\
K>0.$ If the initial datum $y(0)=y_{0}<-\sqrt{\frac{K}{C}},$ then
the solution to (2.4) goes to $-\infty$ before $t$ tend to
$\frac{1}{-Cy_{0}+\frac{K}{y_{0}}}.$

\section{Blow-up phenomena}
\noindent

In this section, we discuss the blow-up phenomena of the system
(2.1). Firstly, we prove that there exist strong solutions to (2.1)
which do not exist globally in time.\\
\newline
\textbf{Theorem 3.1.} Let $(u_{0}, \rho_{0}) \in H^s\times H^{s-1},
s\geq 2,$ and T be the maximal time of the solution $(u,\rho)$ to
(2.1) with the initial data $(u_{0}, \rho_{0}).$ If there is some
$x_{0} \in \mathbb{S}$ such that $\rho_{0}(x_{0})=0$ and
$$u_{0}^{\prime}(x_{0})=\inf\limits_{x\in
\mathbb{S}}u_{0}^{\prime}(x)<-\sqrt{\frac{e+1}{2(e-1)}E_{0}+|\gamma-A|\sqrt{\frac{8(e+1)}{e-1}}E_{0}^{\frac{1}{2}}},$$
then the corresponding solution to (2.1) blows up in finite time.\\
\newline
\textbf{Proof.} Applying Lemma 2.1 and a simple density argument, we
only need to show that the above theorem holds for some $s\geq 2$.
Here we assume $s=3$ to prove the above theorem.\

Define now $$ m(t):=\inf_{x\in \mathbb{S}}[u_{x}(t,x)], \ \ t\in
[0,T).$$ By Lemma 2.7, we let $\xi(t)\in \mathbb{S}$ be a point
where this infimum is attained. It follows that
$$m(t)=u_{x}(t,\xi(t)) \ \ \text{and} \ \ u_{xx}(t,\xi(t))=0.$$

Differentiating the first equation in (2.2) with respect to $x$ and
using the identity $\partial_x^2G\ast f=G\ast f-f$, we have
\begin{equation}
u_{tx}+(u-\gamma)u_{xx}=-\frac{1}{2}u_{x}^{2}+\frac{1}{2}\rho^{2}+u^{2}
+(\gamma-A)u-G*(u^{2}+\frac{1}{2}u_{x}^{2}+\frac{1}{2}\rho^{2}+(\gamma-A)u).
\end{equation}

Since the map $q(t,\cdot)$ given by (2.3) is an increasing
diffeomorphism of $\mathbb{R},$ there exists a $x(t)\in \mathbb{S}$
such that $q(t,x(t))=\xi(t).$ In particular, $x(0)=\xi(0).$ Note
that $u_{0}^{\prime}(x_{0})=\inf\limits_{x\in
\mathbb{S}}u_{0}^{\prime}(x),$ we can choose $x_{0}=\xi(0).$ It
follows that $x(0)=\xi(0)=x_{0}.$ By Lemma 2.3 and the condition
$\rho_{0}(x_{0})=0,$  we have
$$\rho(t,\xi(t))q_{x}(t,x)=\rho(t,q(t,x(t)))q_{x}(t,x)=\rho_{0}(x(0))=\rho_{0}(x_{0})=0.$$
Thus $\rho(t,\xi(t))=0.$\

Valuating (3.1) at $(t,\xi(t))$ and using Lemma 2.7, we obtain
\begin{equation}
\frac{dm(t)}{dt}\leq
-\frac{1}{2}m^{2}(t)+\frac{1}{2}u^{2}+(\gamma-A)u-(\gamma-A)G*u,
\end{equation}
here we used the relations $G*(u^{2}+\frac{1}{2}u_{x}^{2})\geq
\frac{1}{2}u^{2}$ and $G*\rho^{2}\geq 0.$ Note that
$\|G\|_{L^{1}}=1.$ By Lemma 2.4 and Lemma 2.8, we get
$$\|u\|_{L^{\infty}}^{2}\leq \frac{e+1}{2(e-1)}\|u\|_{H^{1}}^{2}\leq
\frac{e+1}{2(e-1)}E_{0},$$
$$|(\gamma-A)u|\leq |\gamma-A|\|u\|_{L^{\infty}}\leq |\gamma-A|
\sqrt{\frac{e+1}{2(e-1)}}E_{0}^{\frac{1}{2}}$$ and
$$|(\gamma-A)G*u|\leq |\gamma-A|\|G\|_{L^{1}}\|u\|_{L^{\infty}}\leq |\gamma-A|
\sqrt{\frac{e+1}{2(e-1)}}E_{0}^{\frac{1}{2}}.$$ It follows that
\begin{equation}
\frac{dm(t)}{dt}\leq -\frac{1}{2}m^{2}(t)+K,
\end{equation}
where $K=\frac{e+1}{4(e-1)}E_{0}+2|\gamma-A|
\sqrt{\frac{e+1}{2(e-1)}}E_{0}^{\frac{1}{2}}.$ Since
$m(0)<-\sqrt{2K},$ Lemma 2.10 implies
$$\lim_{t\rightarrow T}m(t)=-\infty \ \ \ \text{with} \ \
T=\frac{2u_{0}^{\prime}(x_{0})}{2K-(u_{0}^{\prime}(x_{0}))^{2}}.$$
Applying Lemma 2.6, the solution blows up in finite time. $\hfill{}
\Box$ \\
\newline
\textbf{Theorem 3.2.} Let $(u_{0}, \rho_{0}) \in H^s\times H^{s-1},
s\geq 2,$ and T be the maximal time of the solution $(u,\rho)$ to
(2.1) with the initial data $(u_{0}, \rho_{0}).$ Assume that
$\int_{\mathbb{S}}u_{0}(x)dx=\frac{a_{0}}{2}.$ If there is some
$x_{0} \in \mathbb{S}$ such that $\rho_{0}(x_{0})=0$ and for any
$\epsilon >0,$
$$u_{0}^{\prime}(x_{0})=\inf\limits_{x\in
\mathbb{S}}u_{0}^{\prime}(x)<-\sqrt{\frac{\epsilon+2}{24}E_{0}
+\frac{\epsilon+2}{4\epsilon}a_{0}^{2}+|\gamma-A|\sqrt{\frac{2(\epsilon+2)}{3}E_{0}+\frac{4(\epsilon+2)}{\epsilon}a_{0}^{2}}},$$
then the corresponding solution to (2.1) blows up in finite time.\\
\newline
\textbf{Proof.} By Lemma 2.5, we have
$\int_{\mathbb{S}}u(t,x)dx=\frac{a_{0}}{2}.$ Using Lemma 2.4 and
Lemma 2.9, we obtain
$$\|u\|_{L^{\infty}}^{2}\leq \frac{\epsilon+2}{24}E_{0}+\frac{\epsilon+2}{4\epsilon}a_{0}^{2},$$
$$|(\gamma-A)u|\leq |\gamma-A|\|u\|_{L^{\infty}}\leq |\gamma-A|\sqrt{\frac{\epsilon+2}{24}E_{0}+\frac{\epsilon+2}{4\epsilon}a_{0}^{2}}$$ and
$$|(\gamma-A)G*u|\leq |\gamma-A|\|G\|_{L^{1}}\|u\|_{L^{\infty}}\leq |\gamma-A|\sqrt{\frac{\epsilon+2}{24}E_{0}+\frac{\epsilon+2}{4\epsilon}a_{0}^{2}}.$$
Following the similar proof in Theorem 3.1, we have
\begin{equation}
\frac{dm(t)}{dt}\leq -\frac{1}{2}m^{2}(t)+K,
\end{equation}
where
$K=\frac{\epsilon+2}{48}E_{0}+\frac{\epsilon+2}{8\epsilon}a_{0}^{2}+|\gamma-A|\sqrt{\frac{\epsilon+2}{6}E_{0}+\frac{\epsilon+2}{\epsilon}a_{0}^{2}}.$
Following the same argument as in Theorem 3.1, we deduce that the
solution blows up in finite time. $\hfill{}
\Box$ \\
\

Letting $a_{0}=0$ and $\epsilon\rightarrow 0$ in Theorem 3.2, we
have the following result.\\
\newline
\textbf{Corollary 3.1.} Let $(u_{0}, \rho_{0}) \in H^s\times
H^{s-1}, s\geq 2,$ and T be the maximal time of the solution
$(u,\rho)$ to (2.1) with the initial data $(u_{0}, \rho_{0}).$
Assume that $\int_{\mathbb{S}}u_{0}(x)dx=0.$ If there is some $x_{0}
\in \mathbb{S}$ such that $\rho_{0}(x_{0})=0$ and
$$u_{0}^{\prime}(x_{0})=\inf\limits_{x\in
\mathbb{S}}u_{0}^{\prime}(x)<-\sqrt{\frac{E_{0}}{12}+2|\gamma-A|\sqrt{\frac{E_{0}}{3}}},$$
then the corresponding solution to (2.1) blows up in finite time.\\
\newline
\textbf{Remark 3.1.} Note that the system (2.1) is variational under
the transformation $(u,x)\rightarrow(-u,-x)$ and $(\rho,
x)\rightarrow(\rho, -x)$ even $\gamma=0.$ Thus, we can not get a
blow up result according to the parity of the initial data $(u_{0},
\rho_{0})$ as we usually do.\\
\

Next, we will give more insight into the blow-up mechanism for the
wave-breaking solution to the system (2.1), that is the blow-up rate
for strong solutions to (2.1).\\
\newline
\textbf{Theorem 3.3.} Let $(u,\rho)$ be the solution to the system
(2.1) with the initial data $(u_{0}, \rho_{0}) \in H^s\times
H^{s-1}, s\geq 2,$ satisfying the assumption of Theorem 3.1, and T
be the maximal time of the solution $(u,\rho).$ Then, we have
$$ \lim_{t\rightarrow T}(T-t)\inf_{x\in \mathbb{S}}u_{x}(t,x)=-2.$$
\newline
\textbf{Proof.} As mentioned earlier, here we only need to show that
the above theorem holds for $s=3$. \

Define now $$ m(t):=\inf_{x\in \mathbb{S}}[u_{x}(t,x)], \ \ t\in
[0,T).$$ \

By the proof of Theorem 3.1, we have there exists a positive
constant $K=K(E_{0}, \gamma, A)$ such that
\begin{equation}
-K\leq\frac{d}{dt}m+\frac{1}{2}m^{2}\leq K \ \ \ a.e.\ \ on\ (0,T).
\end{equation}
Let $\varepsilon\in(0,\frac{1}{2})$. Since
$\liminf\limits_{t\rightarrow T}m(t)= -\infty$ by Theorem 3.1, there
is some $t_{0}\in (0,T)$ with $m(t_{0})<0$ and
$m^{2}(t_{0})>\frac{K}{\varepsilon}$. Since $m$ is locally
Lipschitz, it is then inferred from (3.5) that
\begin{equation}
m^{2}(t)>\frac{K}{\varepsilon}, \ \ \ \ t\in [t_{0}, T).
\end{equation}\

A combination of (3.5) and (3.6) enables us to infer
\begin{equation}
\frac{1}{2}+\varepsilon\geq-\frac{\frac{dm}{dt}}{m^{2}}\geq\frac{1}{2}-\varepsilon
\ \ \ a.e.\ on \ (0,T).
\end{equation}
Since $m$ is locally Lipschitz on $[0,T)$ and (3.6) holds, it is
easy to check that $\frac{1}{m}$ is locally Lipschitz on
$(t_{0},T).$ Differentiating the relation $m(t)\cdot
\frac{1}{m(t)}=1,\ t\in (t_{0}, T),$ we get
$$\frac{d}{dt}(\frac{1}{m})=-\frac{\frac{dm}{dt}}{m^{2}} \ a.e.\ on\
(t_{0}, T),$$ with $\frac{1}{m}$ absolutely continuous on
$(t_{0},T).$ For $t\in (t_{0},T)$. Integrating (3.7) on $(t,T)$ to
obtain $$
(\frac{1}{2}+\varepsilon)(T-t)\geq-\frac{1}{m(t)}\geq(\frac{1}{2}-\varepsilon)(T-t),\
t\in (t_{0}, T),$$ that is,$$
\frac{1}{\frac{1}{2}+\varepsilon}\leq-m(t)(T-t)\leq\frac{1}{\frac{1}{2}-\varepsilon},\
t\in(t_{0},T).$$ By the arbitrariness of
$\varepsilon\in(0,\frac{1}{2})$ the statement of Theorem 3.3
follows. $\hfill{} \Box$

\section{Global Existence}
\noindent

In this section, we will present a global existence result.\\
\newline
\textbf{Theorem 4.1.} Let $(u_{0}, \rho_{0}) \in H^s\times H^{s-1},
s\geq 2,$ and T be the maximal time of the solution $(u,\rho)$ to
(2.1) with the initial data $(u_{0}, \rho_{0}).$ If $\rho_{0}(x)\neq
0$ for all $x\in\mathbb{S}$, then the corresponding solution
$(u,\rho)$
exists globally in time.\\
\newline
\textbf{Proof.} Define now $$ m(t):=\inf_{x\in
\mathbb{S}}[u_{x}(t,x)], \ \ t\in [0,T).$$ By Lemma 2.7, we let
$\xi(t)\in \mathbb{S}$ be a point where this infimum is attained. It
follows that
$$m(t)=u_{x}(t,\xi(t)) \ \ \text{and} \ \ u_{xx}(t,\xi(t))=0.$$ Since the map $q(t,\cdot)$ given by (2.3) is an increasing
diffeomorphism of $\mathbb{R},$ there exists a $x(t)\in \mathbb{S}$
such that $q(t,x(t))=\xi(t).$\

Set $m(t)=u_{x}(t,\xi(t))=u_{x}(t,q(t,x(t)))$ and
$\alpha(t)=\rho(t,\xi(t))=\rho(t,q(t,x(t))).$ Valuating (3.1) at
$(t,\xi(t))$ and using Lemma 2.7, we obtain
\begin{equation}
m^{\prime}(t)=-\frac{1}{2}m^{2}(t)+\frac{1}{2}\alpha^{2}(t)+f \ \
\text{and} \ \ \alpha^{\prime}(t)=-m(t)\alpha(t),
\end{equation}
where $f=u^{2}
+(\gamma-A)u-G*(u^{2}+\frac{1}{2}u_{x}^{2}+\frac{1}{2}\rho^{2}+(\gamma-A)u).$
By Lemma 2.4, Lemma 2.8 and $\frac{1}{2\sinh\frac{1}{2}}\leq G(x)
\leq \frac{\cosh \frac{1}{2}}{2\sinh\frac{1}{2}},$ we have
\begin{eqnarray*}
|f|&\leq&
\|u\|_{L^{\infty}}^{2}+2|\gamma-A|\|u\|_{L^{\infty}}+\|G\|_{L^{\infty}}\|u^{2}+\frac{1}{2}u_{x}^{2}+\frac{1}{2}\rho^{2}\|_{L^{1}}\\
&\leq&
\frac{e+1}{2(e-1)}E_{0}+2|\gamma-A|\sqrt{\frac{e+1}{2(e-1)}}E_{0}^{\frac{1}{2}}+\frac{\cosh
\frac{1}{2}}{2\sinh\frac{1}{2}}E_{0}:=c_{1}
\end{eqnarray*}\

By Lemmas 2.2-2.3, we know that $\alpha(t)$ has the same sign with
$\alpha(0)=\rho_{0}(x_{0})$ for every $x\in \mathbb{R}$. Moreover,
there is a constant $\beta>0$ such that
$|\alpha(0)|=\inf\limits_{x\in\mathbb{S}}|\rho_{0}(x)|\geq\beta>0$
because of $\rho_{0}(x)\neq 0$ for all $x\in\mathbb{S}.$ Next, we
consider the following Lyapunov positive function
\begin{equation}
w(t)=\alpha(0)\alpha(t)+\frac{\alpha(0)}{\alpha(t)}(1+m^{2}(t)), \ \
\ t\in[0,T).
\end{equation}
Letting $t=0$ in (4.2), we have
$$w(0)\leq
\|\rho_{0}\|_{L^{\infty}}^{2}+1+\|u_{0}^{\prime}(x)\|_{L^{\infty}}^{2}:=c_{2}.$$
Differentiating (4.2) with respect to $t$ and using (4.1), we obtain
\begin{eqnarray*}
w^{\prime}(t)&=&\frac{\alpha(0)}{\alpha(t)}\cdot
2m(t)(f+\frac{1}{2})\\
&\leq& \frac{\alpha(0)}{\alpha(t)}(1+m^{2}(t))(|f|+\frac{1}{2})\\
&\leq& w(t)(c_{1}+\frac{1}{2}).
\end{eqnarray*}
By Gronwall's inequality, we have $$w(t)\leq
w(0)e^{(c_{1}+\frac{1}{2})t}\leq c_{2}e^{(c_{1}+\frac{1}{2})t}$$ for
all $t\in [0,T).$ On the other hand,
$$w(t)\geq 2\sqrt{\alpha^{2}(0)(1+m^{2}(t))}\geq 2 \beta|m(t)|, \ \
\ \forall \ \ t\in [0,T).$$ Thus,
$$
|m(t)|\leq \frac{1}{2\beta}w(t)\leq
\frac{c_{2}}{2\beta}e^{(c_{1}+\frac{1}{2})t}
$$ for all $t \in [0,T)$. It follows that
$$\liminf_{t\rightarrow T}m(t)\geq
-\frac{c_{2}}{2\beta}e^{(c_{1}+\frac{1}{2})T}.$$ This completes the
proof by using Lemma 2.6.


\begin{thebibliography}{99}
\small

\bibitem{zhu} M. Zhu, J. Xu, On the wave-breaking phenomena for the periodic two-component
Dullin-Gottwald-Holm system, J. Math. Anal. Appl. 391 (2012)
415-428.

\bibitem{Ivanov}R. Ivanov, Two-component integrable systems modelling shallow water waves: the constant vorticity case,
Wave Motion 46 (2009) 389-396.

\bibitem{DGH} H. R. Dullin, G. A. Gottwald, D. D. Holm, An integral shallow water equation with linear and nonlinear dispersion,
 Phys. Rev. Lett. 87 (2001) 4501-4504.

\bibitem{yin}Z. Yin, Well-posedness, blowup, and global existence for an integrable shallow water equation,
Discrete Contin. Dyn. Syst. 11 (2004) 393-411.

\bibitem{shen}C. Shen, L. Tian, A. Gao, Optimal control of the viscous
Dullin-Gottwalld-Holm equation, Nonlinear Analysis: Real World
Applications 11 (2010) 480-491.

\bibitem{ai}X. Ai, G. Gui, On the inverse scattering problem and the low
regularity solutions for the Dullin-Gottwald-Holm equation,
Nonlinear Anal. Real World Appl. 11 (2010) 888-894.

\bibitem{guo}Z. Guo, L. Ni, Wave breaking for the periodic weakly
dissipative Dullin-Gottwald-Holm equation, Nonlinear Anal. 74 (2011)
965-973.

\bibitem{meng}Q. Meng, B. He, Y. Long, Z. Li, New exact periodic wave
solutions for the Dullin-Gottwald-Holm equation, Appl. Math. Comput.
218 (2011) 4533-4537.

\bibitem{sun}B. Sun, Maximum principle for optimal distributed control
of the viscous Dullin-Gottwald-Holm equation, Nonlinear Anal. Real
World Appl. 13 (2012) 325-332.

\bibitem{liu}Y. Liu, Global existence and blow-up solutions for a
nonlinear shallow water equation, Math. Ann. 335 (2006) 717-735.

\bibitem{tian}L. Tian, G. Gui, Y. Liu, On the Cauchy problem and the
scattering problem for the Dullin-Gottwald-Holm equation, Comm.
Math. Phys. 257 (2005) 667-701.

\bibitem{yin2}Z. Yin, Global existence and blow-up for a periodic
integrable shallow water equation with linear and nonlinear
dispersion, Dyn. Contin. Discrete Impuls. Syst. Ser. A Math. Anal.
12 (2005) 87-101.

\bibitem{zhang}S. Zhang, Z. Yin, On the blow-up phenomena of the periodic
Dullin-Gottwald-Holm equation, J. Math. Phys. 49 (2008) 1-16.

\bibitem{zhang1}S. Zhang, Z. Yin, Global weak solutions for the
Dullin-Gottwald-Holm equation, Nonlinear Anal. 72 (2010) 1690-1700.

\bibitem{m}O. G. Mustafa, Global conservative solutions of the
Dullin-Gottwald-Holm equation, Discrete Contin. Dyn. Syst. 19 (2007)
575-594.

\bibitem{yan}K. Yan, Z. Yin, On the solutions of the Dullin-Gottwald-Holm equation in Besov
spaces, Nonlinear Anal. Real World Appl. 13 (2012) 2580-2592.

\bibitem{chen1} R. M. Chen and Y. Liu, Wave-breaking and global existence for a
generalized two-component Camassa-Holm system, Int. Math. Res. Not.
(2010), in press.

\bibitem{chen2} R. M. Chen, Y. Liu and Z. Qiao, Stability of solitary waves of a
generalized two-component Camassa-Holm system Comm. Partial
Differential Equations (2010), in press.

\bibitem{gui1} G. Gui and Y. Liu, On the Cauchy problem for the two-component
Camassa-Holm system, Math. Z. 268 (2010) 45-66.

\bibitem{gui2} G. Gui and Y. Liu, On the global existence and wave-breaking
criteria for the two-component Camassa-Holm system, J. Funct. Anal.
258 (2010) 4251-4278.

\bibitem{con}A. Constantin, R. Ivanov, On the integrable two-component
Camassa-Holm shallow water system, Phys. Lett. A 372 (2008)
7129-7132.

\bibitem{escher}J. Escher, O. Lechtenfeld, Z.Y. Yin, Well-posedness and blow-up
phenomena for the 2-component Camassa-Holm equation, Discrete
Contin. Dyn. Syst. 19 (2007) 493-513.

\bibitem{fu}Y. Fu, Y. Liu, C.Z. Qu, Well-posedness and blow-up solution for a
modified two-component periodic Camassa-Holm system with peakons,
Math. Ann. 348 (2010) 415-448.

\bibitem{zhang2}P.Z. Zhang, Y. Liu, Stability of solitary waves and wave-breaking
phenomena for the two-component Camassa-Holm system, Int. Math. Res.
Not. IMRN 11 (2010) 1981-2021.

\bibitem{hu} Q. Hu, On a periodic 2-component Camassa-Holm equation with
vorticity, J. Nonlinear Math. Phys. 18 (2011) 541-556.

\bibitem{Constantin 4}
A. Constantin and J. Escher, Wave breaking for nonlinear nonlocal
shallow water equations, Acta. Math., 181 (1998), 229-243.

\bibitem{Y2}
Z. Yin, On the blow-up of solutions of the periodic Camassa-Holm
equation, Dyn. Contin. Discrete Impuls. Syst. Ser. A Math. Anal. 12
(2005), 375-381.

\bibitem{hu1} Q. Hu and Z. Yin, Blowup phenomena for a new periodic nonlinearly dispersive wave equation,
Math. Nachr. 283 (11) (2010) 1613-1628.

\bibitem{zhou} Y. Zhou, Blow-up of solutions to a nonlinear dispersive rod equation, Calc. Var. Partial Differential Equations 25 (2005)
63-77.


\end{thebibliography}
\end{document}